\documentclass[12pt]{article}
\usepackage{latexsym}
\usepackage{amsmath, amsthm}
\usepackage{amsfonts}
\usepackage{mathrsfs}
\newtheorem{prop}{Proposition}[section]
\newtheorem{lem}{Lemma}[section]
\newtheorem{thm}{Theorem}[section]

\newcommand{\diag}{\mathrm{diag}\,}

\newcommand{\Rs}{\mathbb{R}}

\newcommand{\Lc}{{\cal L} }

\newcommand{\G}{\mbox{gal}}
\newcommand{\bz}{{\bf 0} }

\newcommand{\beq}{ \begin{equation} }
\newcommand{\eeq}{ \end{equation} }
\newcommand{\bt}{ \begin{tabular} }
\newcommand{\et}{ \end{tabular} }

\begin{document}

\bibliographystyle{plain}
\title{A Note on Global Positioning System (GPS) and Euclidean Distance Matrices }
\vspace{0.3in}
        \author{ A. Y. Alfakih
  \thanks{E-mail: alfakih@uwindsor.ca}
  \\
          Department of Mathematics and Statistics \\
          University of Windsor \\
          Windsor, Ontario N9B 3P4 \\
          Canada
}

\date{\today}
\maketitle


\vspace{0.1in}

\begin{abstract}

Let $D$ be an $n \times n$ Euclidean distance matrix with embedding dimension $r$;
and let $d^m \in \Rs^n$ be a given vector. 
In this note, we consider the problem of finding a vector $y \in \Rs^n$,
that is closest to $d^m$ in Euclidean norm, such such that
the augmented matrix
$\left[ \begin{array}{cc} 0 & y^T \\ y & D \end{array} \right]$ is itself an EDM 
with embedding dimension $r$. This problem is motivated by applications in global 
positioning system (GPS). We present a fault detection criterion and 
three algorithms: one for the case $n=4$, and two for the case $n \geq 5$.

\end{abstract}

\section{Introduction}

Recently, several publications \cite{bea24, kg23, ing25} employed   
Euclidean distance matrices (EDMs) to address various problems related to   
the Global Positioning System (GPS). This note follows in the spirit of
these works, aiming to highlight the potential usefulness of EDMs theory in  
the mathematics of GPS. 
EDMs have found applications across a range of fields, including  
molecular conformation 
theory \cite{ch88}, the statistical theory of multidimensional scaling \cite{bg97}, 
wireless sensor networks \cite{kw10}, and the rigidity theory of 
bar-and-joint frameworks \cite{alf18m}.

The Global Positioning System (GPS) is a satellite-based navigation system
that allows users to determine their position anywhere on Earth. 
It consists of a constellation of at least  
24 operational satellites distributed across six orbital planes, at an  
altitude of  about 20,200 km above Earth's surface. 
These satellites are arranged to ensure that at least four are visible from
any point on Earth's surface at all times.

Let $\rho^m_i$ denote the {\em pseudorange} to satellite $i$, i.e.,  
the measured distance between the receiver and satellite $i$. The pseudorange 
$\rho^m_i$ differs from the true geometric range $\rho_i$ 
due to factors such as clock
synchronization error, atmospheric effects, relativistic effects 
and other sources of error (e.g., multipath, receiver noise, etc.). 
Relativistic and atmospheric errors can be effectively corrected using
established models and techniques. 
However, clock error arising from unsynchronized satellite
and receiver clocks, and random error such as those caused by receiver noise
or multipath effects, are more difficult to correct.
Due to their unpredictable nature, these errors can significantly 
degrade positioning accuracy. 
Effectively mitigating them is essential for enhancing the overall accuracy and 
reliability of GPS-based systems. Thus, we can assume that 

\[
(\rho^m_i)^2 = (\rho_i)^2 + \epsilon_i,
\] 
where $\epsilon_i$ is the error associated with satellite $i$ assumed to be 
orders of magnitude smaller than $(\rho_i)^2$. 

This motivates the following problem. 
Let $p^1,\ldots,p^n \in \Rs^r$, with $n \geq r+1$,  
denote the known positions of $n$ satellites, 
\footnote{Although in the actual GPS problem the dimension is $r=3$, we prefer
to maintain generality by using arbitrary $r$.} 
and assume that $p^1,\ldots,p^n$ affinely span $\Rs^r$. 
Let $d^m = (d^m_i = (\rho^m_i)^2 ) \in \Rs^n$, be a given vector, 
and define the $n \times n$ matrix $D=(d_{ij})$ by 
\beq \label{defD} 
 d_{ij} = || p^i - p^j||^2, 
\eeq
where $||.||$ denotes the Euclidean norm. 
Thus $D$ is an EDM of embedding dimension $r$.  
The goal is to find the vector $y \in \Rs^n$ closest to $d^m$ in Euclidean norm,
and then to determine the receiver's position $q \in \Rs^r$, 
such that $y_i = || q - p^i||^2$ for $i=1,\ldots,n$. 
In other words, we seek to solve the optimization problem 
\beq \label{prob1} 
 \begin{tabular}{lll} 
 $\underset{y}{\min}$ & $||y - d^m||^2$ &  \\  
subject to &  $\left[ \begin{array}{cc} 
     0 & y^T \\ y  & D \end{array} \right]$ & is an EDM of embedding dimension $r$. 
\end{tabular}
\eeq 
Once the optimal solution $y^*$ is found, the receiver's position $q$ 
can be recovered via a simple formula. 
We present three algorithms for solving (\ref{prob1}): one for the case $n=4$, and 
two for the case $n \geq 5$. In addition, we provide  a necessary and sufficient
condition for $d^m$ to be self-consistent, meaning that  
$y^* =d^m$ is already the optimal solution of problem (\ref{prob1}).

 It is worth noting that the receiver's position $q$, and hence $y$, can also 
 be obtained by solving the following unconstrained nonlinear optimization
problem
\beq \label{nlp} 
\min_q \sum_{i=1}^n (||p^i-q||^2 - d^m_i)^2. 
\eeq 
As we will show, for $n \geq 5$, the optimality condition of 
problem (\ref{prob1}) leads to an unconstrained minimization problem
equivalent to problem (\ref{nlp}).

\subsection{Notation}

We summarize, here, the notation used throughout this note.
$e$ denotes the vector of all 1's of the appropriate dimension. 
$I$ is the identity matrix of the appropriate dimension. 
$\bz$ denotes the zero vector or the zero matrix of appropriate dimension.
 For a symmetric matrix $A$, the notation $A \succeq \bz$ means that $A$ is positive 
semidefinite, and $A^{\dag}$ is the Moore-Penrose inverse of $A$.  
Finally, $||x|| = \sqrt{x^Tx}$ denotes the Euclidean norm of $x$. 

\section{Euclidean Distance Matrices}

In this section, we present the results of the theory of EDMs most relevant
to this note. 
For a comprehensive discussion, refer to the monograph \cite{alf18m}.

An $n \times n$ matrix $D$ is called a {\em Euclidean distance matrix (EDM)}
if there exist points $p^1,\ldots,p^n$ in some Euclidean space such that
\[
d_{ij} = ||p^i - p^j||^2 \mbox{ for all } i,j=1,\ldots,n.
\] 
These points $p^1,\ldots,p^n$ are called the {\em generating points} of $D$, and 
the dimension of their affine span is called the {\em embedding dimension}
of $D$.  Let the embedding dimension of an EDM $D$ be $r$. 
We assume that $p^1,\ldots,p^n$, the generating points of $D$, 
lie in $\Rs^r$ and affinely span $\Rs^r$. 
Define the $n \times r$ matrix 
\beq \label{defP} 
P = \left[ \begin{array}{c} (p^1)^T \\ \vdots \\ (p^n)^T \end{array} \right],
\eeq
where each $p^i \in \Rs^r$  is a row of $P$. 
$P$ is called a {\em configuration matrix} of $D$. Since the points 
$p^1,\ldots, p^n$ affinely span $\Rs^r$, $P$ has full column rank; i.e.,
rank $P = r$. Without loss of generality, we assume
that the origin coincides with the centroid of $p^1,\ldots,p^n$.   
Thus, 
\beq \label{pe}
P^T e = \bz,
\eeq 
where $e \in \Rs^n$ is the vector of all 1's. 
Let $V$ be an $n \times (n-1)$ matrix such that 
\beq \label{defVQ}  
Q=[e/\sqrt{n} \;\;\;\;  V]
\eeq 
is an $n \times n$ orthogonal matrix. 

Clearly, an EDM $D$ is symmetric with zero diagonal, and nonnegative
offdiagonal entries.  
Among the various characterizations of EDMs, the following two are most relevant
to our purposes. 

\begin{thm} \label{EDMthmV} \cite{gow85,sch35,yh38} 
Let $D$ be an $n \times n$ symmetric 
matrix with zero diagonal. Then $D$ is an EDM if and only if $(-V^T D V)$ is a 
positive semidefinite matrix. Moreover, the embedding dimension of $D$ is equal 
to rank $(-V^TDV)$. 
\end{thm}
The $(n-1) \times (n-1)$ matrix $X = - \frac{1}{2} V^T DV$ is called the 
{\em projected Gram matrix} of $D$. 
 
\begin{thm} \label{EDMthmX} \cite{gow85,sch35,yh38} 
Let $D$ be an $n \times n$ symmetric 
matrix with zero diagonal and let $d \in \Rs^n$. Then, the augmented matrix 
$\left[ \begin{array}{cc} 0 & d^T \\ d & D \end{array} \right]$ is an
EDM with embedding dimension $r$ if and only if 
\[
d e^T + e d^T - D 
\]
is a positive semidefinite matrix with rank $r$. 
\end{thm}

The Gram matrix of the points $p^1, \ldots, p^n$ is given by  
\beq \label{defB} 
B= P P^T. 
\eeq   
Note that our assumption $P^T e = \bz$ implies that $B e = \bz$. 
Matrix $B$ is positive semidefinite with rank $r$. 
The Gram matrix $B$ and the projected Gram matrix $X$ are related by  
\[
X = V^T B V \mbox{ and } B= V X V^T.
\]
Let $B^{\dag}$ denote the Moore-Penrose inverse of $B$. Then it is
easy to verify that 
\beq \label{defBi}
B^{\dag} = V X^{\dag} V^T = P (P^TP)^{-2}P^T. 
\eeq

Let $\diag(B)$ denote the vector formed from the diagonal entries of $B$. Then 
the Gram matrix $B$ and its associated EDM $D$ are related by  
\beq \label{DB} 
D = \diag(B) \; e^T + e \;  (\diag(B))^T - 2 B, 
\eeq
and 
\beq \label{BD}
B = - \frac{1}{2} J D J,
\eeq
where $J = VV^T = I - ee^T/n$. 
Let $P$ be a configuration matrix of an EDM $D$ with embedding dimension $r$.
Then the {\em Gale space} of $D$ is defined as
\beq \label{defZ}
\G(D) = \mbox{null space of } 
\left[ \begin{array}{c} P^T \\ e^T  \end{array} \right]. 
\eeq 
Furthermore, any $n \times (n-r-1)$ matrix $Z$ whose columns form a basis of
$\G(D)$ is called a {\em Gale matrix} of $D$.
It should be pointed out that the columns of $Z$ express the affine dependency
of the points $p^1,\ldots,p^n$.

\begin{lem} \cite[Lemma 3.8]{alf18m} \label{VU}
Let $D$ be an EDM with embedding dimension $r$ and let $X$ be its projected 
Gram matrix. Further, let $U$ be 
the $(n-1) \times (n-1-r)$ matrix whose columns form a basis of the null
space of $X$, and assume that the configuration matrix $P$ of $D$ satisfies
$P^T e = \bz$. Then  $Z=VU$ is a Gale matrix of $D$.  

\end{lem} 

\section{Main Results}

Given an $n \times n$ EDM $D$ with embedding dimension $r$ and a 
vector $y \in \Rs^n$, $y$ is said to be {\em self-consistent} if  
the $(n+1) \times (n+1)$ matrix 
 \beq \label{np1EDM}
\left[ \begin{array}{cc} 0 & y^T \\ y  & D \end{array} \right]  
\eeq 
 is an EDM with embedding dimension $r$. 
Otherwise, $y$ is said to be {\em faulty}. 
By Theorem \ref{EDMthmX}, $y$ is self-consistent if and
only if the $n \times n$ matrix  
\beq \label{matX}
y e^T + e y^T - D
\eeq
is positive semidefinite with rank $r$.

Let $Q$ be the orthogonal matrix defined in (\ref{defVQ}). 
Multiplying the matrix in (\ref{matX}) from the left by $Q^T$ and from the
right by $Q$, we obtain  
\beq \label{mat1} 
\left[ \begin{array}{cc} 
 2 e^Ty - e^TDe/n & 
     \sqrt{n} \; (y^T - e^TD/n) V  \\ 
     \sqrt{n} \; V^T (y - De/n) & -V^T D V \end{array} \right]. 
\eeq 
From equation (\ref{DB}), we have  
$D e/n = (\diag(B) + (e^T \diag(B)/n) \; e)$ and $e^T D e/n = 2 e^T \diag(B)$.
Additionally, by definition, $(-V^TDV) = 2 X$, where $X$ is the projected
Gram matrix associated with $D$. Substituting these into the matrix 
in (\ref{mat1}), it simplifies to 

\beq \label{mat2} 
\left[ \begin{array}{cc} 
 2 e^T (y -b) & \sqrt{n} \; (y - b)^T V  \\ 
     \sqrt{n} \; V^T (y - b) & 2 X \end{array} \right], 
\eeq 
where $b = \diag(B)$. Thus, $y$ is self-consistent if and only if 
this matrix is positive semidefinite with rank $r$. Furthermore,
the optimization problem in (\ref{prob1}) can be reformulated as:  
\beq \label{prob11} 
 \begin{tabular}{lll} 
 $\underset{y}{\min}$ & $||y - d^m||^2$ &  \\  
subject to &  $\left[ \begin{array}{cc} 
 2 e^T (y -b) & \sqrt{n} \; (y - b)^T V  \\ 
     \sqrt{n} \; V^T (y - b) & 2 X \end{array} \right] \succeq \bz$, with rank $r$. 
\end{tabular}
\eeq 
To simplify notation, and to facilitate solving this problem, we consider 
the cases $n =4$ and $n \geq 5$ separately, always keeping in mind that $r=3$. 
The key distinction is that the Gale space is
trivial for $n =4$ but nontrivial for $n \geq 5$, making the Gale matrix
significant in the latter case. 

We will find it useful to define the function  
\beq \label{defk} 
\kappa_n(y) = \frac{4}{n} e^T(y-b)- (y- b)^T B^{\dag} (y - b),   
\eeq 
where $b= \diag(B)$. This function is central to our analysis in  
this note.  

\subsection{The Case of $n=4$}

For $n=4$, the projected Gram matrix $X$ is $3 \times 3$,  
positive semidefinite with rank $3$; i.e., $X$ is positive definite.
Thus, $X^{\dag} = X^{-1}$ and $B^{\dag} = V X^{-1} V^T$.  
Define the matrix  $M= \left[ \begin{array}{cc} 
  1 & \bz \\ - X^{-1} V^T (y - b)  & I \end{array} \right]$.  
Multiplying the matrix in (\ref{mat2}), after setting $n=4$,
 from the left by $M^T$ and from the
right by $M$, we obtain  
\beq \label{QXQ} 
  \left[ \begin{array}{cc} 
 2 e^T (y -b) - 2 (y-b)^T B^{\dag} (y -b) & \bz  \\ 
     \bz & 2 X \end{array} \right]. 
\eeq 
 
The term $2 e^T (y -b) - 2 (y-b)^T B^{\dag} (y -b)$ is called the
{\em Schur complement} \cite{bv04} of $2 X$ in the matrix in (\ref{mat2}).  
Since $M$ is nonsingular and using the definition
of $\kappa_n(y)$ in (\ref{defk}), it follows that the matrix 
$\left[ \begin{array}{cc} 0 & y^T \\ y & D \end{array} \right]$ is:  
\begin{itemize}
\item an EDM of embedding dim = $3$ iff $\kappa_4(y) = 0$, 
\item an EDM of embedding dim = $4$ iff $\kappa_4(y) > 0$,  
\item  not an EDM iff $\kappa_4 (y) < 0$, 
\end{itemize}
Let $d^m \in \Rs^4$ denote the squared pseudoranges from 4 satellites.

\begin{prop} \label{p1} 
For a given vector $ d^m \in \Rs^4$ and a $4 \times 4$ EDM  $D$ with  
embedding dimension 3,  $d^m$ is self-consistent if and only if
$\kappa_4(d^m) = 0$. i.e., iff  
\begin{eqnarray*}
e^T(d^m -b) & = &  (d^m- b)^T B^{\dag} (d^m - b), \\  
            & = &  (d^m- b)^T P (P^TP)^{-2} P^T (d^m - b),  
\end{eqnarray*}
\end{prop}

As an application of Proposition \ref{p1}, consider the case where the errors
in $d^m$ are constant, i.e., $d^m = d + \delta e$, where 
$d$ is the unknown true geometric range, and $\delta$ is a scalar. 
Since $d=d^m - \delta e$ is self-consistent, 
Proposition \ref{p1} implies that
\[
\delta 
= \frac{1}{4} \kappa_4(d^m) = \frac{1}{4} (e^T(d^m -b)-(d^m- b)^T B^{\dag}( d^m -b)).
\] 

As a result, the optimization problem in (\ref{prob11}) is equivalent to 

\beq \label{prob3} 
 \begin{tabular}{lll} 
 $\underset{y}{\min}$ & $||y - d^m||^2$ &  \\  
subject to &  $ \kappa_4(y)= e^T(y-b) - (y- b)^T B^{\dag} (y - b) = 0$,  
\end{tabular}
\eeq 
This is a quadratic programming problem with a single
quadratic equality constraint, solvable by the method 
described in \cite{hma10,acg11} \footnote{The referenced work \cite{hma10} addresses 
the same GPS problem and arrives at an optimization problem similar to (\ref{prob4})
without employing Euclidean distance matrices (EDMs).}. 
Below, we specialize this method to our case.  

Let $B^{\dag} = S \Lambda S^T$ be the spectral decomposition of $B^{\dag}$,
where $\Lambda$ is the diagonal matrix consisting of the eigenvalues of
$B^{\dag}$; and $S$ is the orthogonal matrix of
the corresponding eigenvectors. Assume that the eigenvalues of $B^{\dag}$
are $\mu_1 \geq \mu_2 \geq  \mu_3 > \mu_4 = 0$. 

Define the transformation  
\beq 
y = Sx + d^m . 
\eeq

Then the optimization problem in (\ref{prob3}) reduces to   
\beq \label{prob4} 
 \begin{tabular}{lll} 
 $\underset{x}{\min}$ & $x^Tx$ &  \\  
subject to &  $ x^T \Lambda x - 2 c^T x - \kappa_4(d^m) = 0$, 
\end{tabular}
\eeq 
where $c= -\Lambda S^T (d^m-b) + S^T e/2$. 
Let $s^i$ denote the $i$th column  of $S$. Then, since $B^{\dag} e =\bz$, we can
set $s_4 = e/2$.    
Hence, 
\[
c= \left[ \begin{array}{c} -\mu_1 \; {s^1}^T(d^m-b) \\ -\mu_2 \; {s^2}^T(d^m-b) \\
                          -\mu_3  \; {s^3}^T(d^m-b) \\ 1 \end{array} \right]. 
\]
To avoid pathological cases, we assume that $c_1 \neq 0$, i.e., $d^m - b$ is 
not orthogonal to $s^1$, the eigenvector of $B^{\dag}$ corresponding to 
its largest eigenvalue $\mu_1$. 

The Lagrangian \cite{bv04} for this problem is  
\[
\Lc (x,\lambda) = x^Tx - \lambda (x^T \Lambda x - 2 c^T x - \kappa_4(d^m) ), 
\]
where $\lambda$ is the Lagrange multiplier. 
The first-order Karush-Kuhn-Tucker (KKT) 
conditions require the gradient of the Lagrangian
with respect to $x$ and $\lambda$ to vanish at a stationary point. Thus  

\begin{eqnarray}
\nabla_x {\Lc} & = &  (I - \lambda \Lambda) x  + \lambda c  = 0 \label{op1} \\
\nabla_{\lambda} \Lc & = &  -x^T \Lambda x + 2 c^T x +\kappa_4(d^m)  = 0. \label{op2} 
\end{eqnarray}
The Hessian of the Lagrangian is 
\beq
\nabla^2 \Lc = \left[ \begin{array}{cc} \nabla_{xx} & \nabla_{ x\lambda} \\
               \nabla_{\lambda x} & \nabla_{\lambda \lambda} \end{array} \right]   
= \left[ \begin{array}{cc} I - \lambda \Lambda & -\Lambda x + c  \\
               - x^T \Lambda + c^T & 0 \end{array} \right].   
\eeq 
The second-order sufficient KKT condition requires the Hessian to be
positive definite on the constraint tangent space, i.e., 
$v^T (I - \lambda \Lambda ) v > 0$ for all 
$v \neq \bz$ such that $(-\Lambda x + c) v = 0$. 
This is satisfied if the optimal Lagrange multiplier $\lambda^*$ 
satisfies 
\beq \label{sos} 
 \lambda^* < \frac{1}{\mu_1} 
\eeq 
ensuring that $I - \lambda^* \Lambda$ is positive definite.  
 The constraint qualification condition  
\[
\Lambda x  - c \neq \bz
\]
holds since $c_4 = 1$, i.e., $c$ does not lie in the column space of $\Lambda$. 

For $\lambda < 1/ \mu_1$, the matrix $I - \lambda \Lambda$ is positive
definite and thus nonsingular.
From the KKT condition  (\ref{op1}), we obtain 

\beq \label{xlam}
x = - \lambda (I - \lambda \Lambda)^{-1} c.  
\eeq
Substituting this into the KKT condition (\ref{op2}) yields   
\[
g(\lambda) =\lambda^2 c^T (I - \lambda \Lambda)^{-1} \Lambda 
(I - \lambda \Lambda)^{-1} c +2 \lambda c^T (I - \lambda \Lambda)^{-1} c 
                                                 - \kappa_4(d^m) = 0, 
\]
or equivalently  
\[
g(\lambda)= \sum_{i=1}^3 
(\lambda^2 \frac{{\mu_i c_i}^2}{(1 - \lambda \mu_i)^2} 
+ 2 \lambda \frac{{c_i}^2}{(1 - \lambda \mu_i)}) +2 \lambda -\kappa_4(d^m) = 0.  
\]
Adding and subtracting the term $\sum_{i=1}^3 \frac{c_i^2}{\mu_i}$ we rewrite  
\beq 
g(\lambda)= \sum_{i=1}^3 
\frac{c_i^2}{\mu_i (1 - \lambda \mu_i)^2} 
+2 \lambda  - h = 0,  
\eeq 
where $h = \kappa_4(d^m) +  \sum_{i=1}^3 \frac{c_i^2 }{\mu_i}$. 
Recall our assumption that $c_1 \neq 0$.
The function $g(\lambda)$ is strictly increasing 
for $\lambda < 1/\mu_1$ with $g(0) = -\kappa_4(d^m)$. Now, if
$\kappa_4(d^m) = 0$, then $\lambda =0$ is the root of $g(\lambda)$ and the
optimal solution of problem (\ref{prob1}) is $y^* = d^m$ indicating that 
$d^m$ is self-consistent. On the other hand, if   
$\kappa_4(d^m) > 0$, then $g(\lambda)$ has a unique root in 
the interval $(0, 1/\mu_1)$.     
Finally, if $\kappa_4(d^m) < 0$, then $g(\frac{1}{2}\kappa_4(d^m)) < 0$ and thus
$g(\lambda)$ has a unique root in the interval $(\frac{1}{2}\kappa_4(d^m), 0)$.

The root $\lambda^*$ of $g(\lambda) = 0$ can be computed using, for instance,
the bisection method provided in {\em julia}'s package {\em roots.jl} 
\cite{beks17}. 
The optimal solution of problem (\ref{prob1}) is then  
\[
y^* = -\lambda^* S( I -\lambda^* \Lambda)^{-1} c + d^m. 
\]

\section{The Case of $n \geq 5$}

Recall that $X$ is the projected Gram matrix associated with $D$, and that
$X$ is an $(n-1) \times (n-1)$ positive semidefinite matrix with rank $r$. 
Let $[W \;\; U]$ be the $(n-1) \times (n-1)$ orthogonal matrix, where
the columns of the $(n-1) \times r$ submatrix $W$ are the 
eigenvectors of $X$ corresponding to its positive eigenvalues, and the  
columns of the $(n-1) \times (n-1-r)$ submatrix $U$ are the 
eigenvectors of $X$ corresponding to its zero eigenvalues. 
Thus $X = W \Delta W^T$, where $\Delta$ is the $r \times r$ diagonal matrix 
formed from the positive eigenvalues of $X$. 
Let 
$Q' =\left[ \begin{array}{ccc} 1 & \bz & \bz \\ \bz & W & U \end{array} \right]$.
Then obviously, $Q'$ is orthogonal. Moreover, by multiplying the matrix in  
(\ref{mat1}) from the left by ${Q'}^T$ and from the right by $Q'$, we obtain  

\beq \label{mat1n} 
\left[ \begin{array}{ccc} 
 2 e^T(y - b)  & \sqrt{n} \; (y - b)^T VW & \sqrt{n} \; (y- b)^T VU \\ 
   \sqrt{n} \; W^T V^T (y - b) & 2 \Delta & \bz \\ 
                   \sqrt{n} \;    U^T V^T (y - b) & \bz & \bz  
 \end{array} \right], 
\eeq 
where we used the fact that $e^TDe/n = 2 e^T b$, where $b=\diag(B)$,
$De/n = b + (e^Tb/n) e$, and since $(-V^TDV) = 2 X$. 
Define the nonsingular matrix 
\[
M'= \left[ \begin{array}{ccc} 
  1 &  \bz & \bz \\  - \frac{\sqrt{n}}{2} \Delta^{-1}W^T V^T (y-b) & I &  \bz  \\
            \bz & \bz &  I \end{array} \right].  
\]
By multiplying the matrix in (\ref{mat1n}) from the left by ${M'}^T$ and from the
right by $M'$, we obtain  
\[
\left[ \begin{array}{ccc} 
2e^T(y-b)-\frac{n}{2} (y - b)^T V W \Delta^{-1}W^T V^T (y - b) &  
  \bz & \sqrt{n} \; (y- b)^T VU \\ 
   \bz  & 2 \Delta & \bz \\ 
               \sqrt{n} \;   U^T V^T (y - b) & \bz & \bz  
 \end{array} \right]
\]
Hence, the matrix in (\ref{mat1n}) is positive semidefinite
with rank $r$ if and only if 
\beq \label{eqs1n} 
\begin{array}{rrl}
  U^T V^T ( y - b) & = & 0 \\
2e^T(y-b)-\frac{n}{2} (y - b)^T V W \Delta^{-1}W^T V^T (y - b) & = & 0. 
\end{array}
\eeq 
Now $W \Delta^{-1} W^T = X^{\dag}$
and $V X^{\dag}V^T = B^{\dag}$. Furthermore, Lemma \ref{VU} implies  
that $VU = Z$ is a Gale matrix.
Hence, (\ref{eqs1n}) is equivalent to 
\begin{eqnarray*} \label{eq1n} 
  Z^T ( y - b) & = & 0 \\
2e^T(y-b)-\frac{n}{2} (y^T -b) B^{\dag} (y - b) & = & 0. 
\end{eqnarray*} 
As a result, $d^m$ is self-consistent if and only if 
\begin{eqnarray*}
  Z^T(d^m - b) & = & 0  \\   
      -\kappa(d^m)= (d^m- b)^T B^{\dag} (d^m - b) -  \frac{4}{n} e^T(d^m -b) & = & 0.
\end{eqnarray*}

Moreover, the optimization problem in (\ref{prob1}) is equivalent to 

\beq \label{prob3n} 
 \begin{tabular}{ll} 
 $\underset{y}{\min}$ & $||y -d^m||^2$  \\  
subject to & $ Z^T(y - b) = 0$  \\   
            &  $ (y- b)^T B^{\dag} (y - b) -  \frac{4}{n} e^T(y -b) = 0$, 
\end{tabular}
\eeq 

Using the definition of Gale matrix in (\ref{defZ}), 
the first constraint implies that
\[
y - b = P x + s e,
\]
for some vector $x \in \Rs^r$ and scalar $s$. Substituting this into the
second constraint and recalling that $B^{\dag} = P (P^TP)^{-2} P^T$, 
we obtain 
\[
(x^T P^T + s e^T ) B^{\dag} (P x + s e) -\frac{4}{n} e^T (P x + s e) =
x^T x - 4 s = 0. 
\]
Now $y-d^m = Px +s e + b-d^m$. Thus

\begin{eqnarray*} 
(y - d^m)^T(y-d^m) &  = &  (Px + s e + b-d^m)^T(Px + s e + b - d^m)  \\  
               & = & x^T P^TP x +2x^T P^T(b-d^m) + n s^2  \\
               & & + 2 s e^T(b-d^m)+ (b-d^m)^T(b-d^m). 
\end{eqnarray*} 

As a result, the optimization problem in (\ref{prob3n}) can be reformulated 
in two ways. First, as the following unconstrained optimization problem 

\beq \label{prob4n} 
 \underset{y}{\min} \;\;\; 
 n (x^Tx)^2/16 + x^T P^TP x + x^Tx \; e^T(b-d^m)/2 +2x^T P^T(b-d^m), 
\eeq
by substituting $s= x^Tx/4$. 
This problem can be solved by any nonlinear optimization solver such as 
the one provided in {\em julia}'s package {\em optim.jl}.
If $x^*$ is the optimal solution of this problem, then the
optimal solution of problem (\ref{prob1}) is 
\[
y^*= P x^* + \frac{1}{4} {x^*}^Tx^* \; e + b.
\]

Second, problem (\ref{prob3n}) can also be formulated as the 
following quadratically constrained quadratic problem 
\beq \label{prob5n} 
 \begin{tabular}{lc} 
 $\underset{x, s }{\min}$ & $x^T P^TP x -2x^T P^T(d^m-b) + n s^2  
                - 2 s e^T(d^m - b)$ \\ 
subject to & $ x^T x - 4 s = 0$  
\end{tabular}
\eeq 
The first-order KKT conditions of this problem are
\begin{eqnarray*}
(P^TP -\lambda I) x -P^T(d^m-b) & = & 0 \\
ns -e^T (d^m - b) +2 \lambda & = & 0, \\
  x^T x - 4 s & = & 0 ,
\end{eqnarray*} 
where $\lambda$ is the Lagrange multiplier. 

Let the  eigenvalues of $P^TP$ be 
$\nu_1 \geq \nu_2 \geq \nu_3 > 0$. 
The second-order sufficient KKT condition requires that matrix  
$ P^TP - \lambda I $ is positive definite, which holds if 
\beq  
\lambda < \nu_3. 
\eeq 
This condition is identical to (\ref{sos}) as the positive eigenvalues 
of $B^{\dag}$ are the reciprocal of those of $P^TP$ \cite{str23}. 

Hence, 
solving for $x$ and $s$ in the above KKT conditions and assuming that
$\lambda < \nu_3$, we obtain 
\beq \label{newx}
x =  ( P^TP - \lambda I)^{-1}  P^T(d^m-b),
\eeq 
\beq \label{news}
s = \frac{1}{n} ( e^T(d^m - b) - 2 \lambda).
\eeq 

Let $P^TP= S' \Lambda' {S'}^T$ be the spectral decomposition of $P^TP$. Then
 substituting (\ref{newx}) and (\ref{news}) into the third KKT condition, 
we obtain 
\[ 
f(\lambda) = (d^m - b)^TP S' (\Lambda' - \lambda I)^{-2}{S'}^T P^T (d^m - b) +
\frac{8}{n} \lambda -\frac{4}{n}e^T(d^m -b) = 0 ,   
\]
or equivalently
\beq
 f(\lambda)= \sum_{i=1}^3 \frac{w_i^2}{(\nu_i - \lambda)^2} 
+ \frac{8}{n} \lambda - h' = 0,   
\eeq 
where $w=S'^T P^T (d^m-b)$ and $h' = 4 e^T (d^m-b)/n$. 
As in the case of $n=4$, to avoid pathological cases we assume that $w_3 \neq 0$.

The function $f(\lambda)$ is strictly increasing 
for $\lambda < \nu_3$. Now  
\begin{eqnarray*} 
 \sum_{i=1}^3 (w_i/\nu_i)^2  & = & (d^m-b)^T P S' (\Lambda')^{-2}{S'}^T P^T(d^m-b) \\ 
     & = & (d^m -b)^T P (P^TP)^{-2} P^T (d^m-b) \\
         & = & (d^m-b)^T B^{\dag} (d^m - b).  
\end{eqnarray*} 
Thus, $f(0) = - \kappa_n(d^m)$. Therefore, if $\kappa_n(d^m) = 0$, then $\lambda=0$
is the root of $f(\lambda)$ and the optimal solution  of (\ref{prob1}) is
$y^*=d^m$ indicating that $d^m$ is self-consistent. On the other hand,  
if $\kappa_n (d^m) > 0$, then $f(\lambda)$ has 
a unique root in the interval $(0, \nu_3)$.     
Finally, if $\kappa_n(d^m) < 0$, then $f(\frac{n}{8} \; \kappa_n(d^m)) < 0$ and thus
$f(\lambda)$ has a unique root in the interval $(\frac{n}{8} \kappa_n(d^m), 0)$.

The root $\lambda^*$ of $f(\lambda) = 0$ can be computed using, for instance,
the bisection method provided in {\em julia}'s package {\em roots.jl}. 
The optimal solution of problem (\ref{prob1}) is then 
\[
y^* = P x^* + s^* e + b,  
\]
where $x^*$ and $s^*$ are given by 
(\ref{newx}) and (\ref{news}) evaluated at $\lambda = \lambda^*$.

\section{Determining the Position of the Receiver $q$}

The position of the receiver $q$ can be determined once the
optimal solution $y^*$ of problem (\ref{prob1}) is obtained. The
$(n+1) \times (n+1)$ Gram matrix for the satellites and the 
receiver is 
\[
\left[ \begin{array}{c} q^T \\ P \end{array} \right] 
\left[ \begin{array}{cc} q &  P^T \end{array} \right]
= \left[ \begin{array}{rr} q^Tq & q^T P^T \\ P q & PP^T \end{array} \right]. 
\]
Hence, using (\ref{DB}), the corresponding EDM is 
\[
 \left[ \begin{array}{rr} 0 & {y^*}^T \\ y^* & D \end{array} \right] =  
\left[ \begin{array}{c} 1 \\ e \end{array} \right] 
\left[ \begin{array}{cc} q^Tq &  b^T \end{array} \right] + 
\left[ \begin{array}{c} q^Tq \\ b \end{array} \right] 
\left[ \begin{array}{cc} 1 &  e^T \end{array} \right]   
- 2  \left[ \begin{array}{rr} q^Tq & q^T P^T \\ P q & PP^T \end{array} \right]. 
\]
From this, we derive
\beq \label{eqq}
2P q = q^Tq \; e + b - y^*. 
\eeq
System of equations (\ref{eqq}) has a solution if and only if
\[
(b-y^*) \mbox{ lies in the column space of } [P \;\; e]. \;\;\;\;\;\;\; (*) 
\]
For $n=4$, condition (*) holds  
trivially since the column space of $[P \;\; e]$ spans all of 
$\Rs^4$. On the other hand, for $n \geq 5$, the definition of the 
Gale matrix $Z$ in (\ref{defZ}) implies that condition (*)  is satisfied
 if and only if $Z^T (b - y^*) = \bz$, which is ensured by 
the first constraint of problem (\ref{prob3n}).  

Multiplying (\ref{eqq})
from the left by $P^T$ and $e^T$, respectively, yields 
\beq \label{feq1} 
q = \frac{1}{2} (P^TP)^{-1} P^T ( b -y^*), 
\eeq 
and 
\beq \label{feq2} 
 q^Tq  = \frac{1}{n} e^T (y^* -b). 
\eeq 
To verify  consistency between (\ref{feq1}) and (\ref{feq2}), note that
$ q^T q = \frac{1}{4} (b-y^*)^T B^{\dag} (b-y^*)$ since 
$P (P^TP)^{-2} P^T = B^{\dag}$. But,  
from the second constraint of problem (\ref{prob3n}), we have  
 $ \kappa_n(y^*)=  \frac{4}{n} e^T(y^* -b)- (y^*- b)^T B^{\dag} (y^* - b) = 0$. 
Thus, $q^T q = \frac{1}{n} e^T (y^* -b)$ confirming that (\ref{feq1}) implies
(\ref{feq2}). 



\begin{thebibliography}{1}

\bibitem{acg11}
C.J.~Albers, F.~Critchley and J.C.~Gower.
\newblock Quadratic minimisation problems in statistics.
\newblock {\em J. Multivariate Analysis.}, 102: 698--713, 2011. 

\bibitem{alf18m}
A.~Y. Alfakih.
\newblock {\em {E}uclidean distance matrices and their applications in rigidity
  theory}.
\newblock Springer, 2018.

\bibitem{bea24}
D.~Beatty.
\newblock Distribution of test statistic for {E}uclidean distance matrices. 
\newblock ar{X}iv 2405.10049v1, 2024.  

\bibitem{bg97}
I.~Borg and P.~Groenen.
\newblock {\em Modern multidimensional scaling, theory and applications}.
\newblock Springer - Verlag , New York, 1997.

\bibitem{beks17}
J.~Bezanson,  A.~Edelman, S.~Karpinski and V.B.~Shah.
\newblock Julia: a fresh approach to numerical computing.
\newblock {\em SIAM review} 59(1), 65--98, 2017. 

\bibitem{bv04}
S.~Boyd and L.~Vandenberghe.  
\newblock {\em Convex optimization}, Cambridge University Press, 2004.

\bibitem{ch88}
G.~M. Crippen and T.~F. Havel.
\newblock {\em Distance Geometry and Molecular Conformation}.
\newblock Wiley, New York, 1988.

\bibitem{gow85}
J.~C. Gower.
\newblock Properties of {E}uclidean and non-{E}uclidean distance matrices.
\newblock {\em Linear Algebra Appl.}, 67:81--97, 1985.

\bibitem{hma10} 
H.~Hmam.
\newblock Quadratic optimisation with one quadratic equality constraint. 
\newblock Electronic Warfare and Radar Division, DSTO-TR--2416, Australia, 2010.  

\bibitem{ing25}
K.~Iiyama, D.~Neamati and G.~Gao.  
\newblock Satellite autonomous clock fault monitoring with inter-satellite       
          ranges using {E}uclidean distance matrices
\newblock ar{X}iv 2505.03820v1. 2025 

\bibitem{kg23}
D.~Knowles and G.~Gao.
\newblock {E}uclidean distance matrix-based rapid fault detection and exclusion. 
\newblock {\em Navigation}, 70 (1), 2023. 

\bibitem{kw10}
N.~Krislock and H.~Wolkowicz.
\newblock Explicit sensor network localization using semidefinite
  representation and facial reduction.
\newblock {\em SIAM J. OPTIM.}, 20:2679--2708, 2010.

\bibitem{sch35}
I.~J. Schoenberg.
\newblock Remarks to {M}aurice {F}r\'{e}chet's article: Sur la d\'{e}finition
  axiomatique d'une classe d'espaces vectoriels distanci\'{e}s applicables
  vectoriellement sur l'espace de {H}ilbert.
\newblock {\em Ann. Math.}, 36:724--732, 1935.

\bibitem{str23}
G.~Strang. 
\newblock {\em Introduction to linear algebra, 6th edition},
\newblock Wellesley-Cambridge Press, Wellesley, MA, 2023.  

\bibitem{yh38}
G.~Young and A.~S. Householder.
\newblock Discussion of a set of points in terms of their mutual distances.
\newblock {\em Psychometrika}, 3:19--22, 1938.

\end{thebibliography}

\end{document}